\documentstyle[11pt,a4wide]{article}

\font \tenmsb=msbm10 scaled \magstep 1
\font \sevenmsb=msbm7 scaled \magstep 1
\font \fivemsb=msbm5 scaled \magstep 1
\newfam \msbfam
\textfont \msbfam=\tenmsb
\scriptfont \msbfam=\sevenmsb
\scriptscriptfont \msbfam=\fivemsb
\def \Bbb#1{\fam \msbfam \relax#1}

\title{\bf  FUNCTION THEORY IN QUANTUM DISC: INTEGRAL REPRESENTATIONS}

\author{\sl S. Sinel'shchikov \and \sl D. Shklyarov
\and \sl L. Vaksman}

\date{\tt Institute for Low Temperature Physics \& Engineering\\
National Academy of Sciences of Ukraine}

\newtheorem{theorem}{Theorem}[section]
\newtheorem{lemma}[theorem]{Lemma}
\newtheorem{proposition}[theorem]{Proposition}
\newtheorem{corollary}[theorem]{Corollary}

\begin{document}

\maketitle
\section{Introduction}

  The theory of von Neumann algebras is a non-commutative analogue of the
function theory of real variable and its far-reaching generalization. This
approach is also applicable for developing a non-commutative analogue of
the function theory of complex variable, as one can see in the classical
work of Arveson \cite{Ar}. However, a progress in non-commutative complex
analysis was inhibited by the absence of substantial examples. The break
was made possible by the fundamental works of L. Faddeev and his
collaborators, V. Drinfeld, M. Jimbo, S. Woronovicz in quantum group theory.
The subject of this paper is the simplest example of a quantum homogeneous
complex manifold, the quantum disc.  Meanwhile, a part of the exposed
results are extensible onto the case of quantum bounded symmetric domains
\cite{SV}.

We obtain non-commutative analogues for some integral representations of
functions in the disc and find a non-commutative analogue of the Plancherel
measure. The Berezin method is applied to produce a formal deformation for
our quantum disc. Our results provide an essential supplement to those of
the well known work \cite{KL} in this field. The initial sections contain
the necessary material on differential and integral calculi in the quantum
disc.

In the cases when the statement in question to be proved via the quantum
group theory is available in the mathematical literature (in particular, in
our electronic preprints \cite{SSV2},\cite{SSV3},\cite{SSV4},\cite{SSV5}), we
restrict ourselves to making exact references.  This approach allows us to
expound all the main results without applying the quantum group theory, and
thus making this text intelligible for a broader class of readers.

The authors are grateful to V. Drinfeld for helpful discussions of a draft
version of this work.

\bigskip

\section{Functions in the quantum disc}

 We assume in the sequel all the vector spaces to be complex, and let $q$
stand for a real number, $0<q<1$.

 The works \cite{NN, KL, F, CK} consider the involutive algebra ${\rm
Pol}({\Bbb C})_q$ given by its generator $z$ and the commutation relation
$$z^*z=q^2zz^*+1-q^2.$$
It arises both in studying the algebras of functions in the quantum disc
\cite{NN, KL, F} and in studying the q-analogues of the Weil algebra
(oscillation algebra) \cite{CK}.

 Any element $f \in{\rm Pol}({\Bbb C})_q$ admits a unique decomposition
$$f=\sum_{jk}a_{jk}(f)z^jz^{*k},\quad a_{jk}\in{\Bbb C}.$$
Moreover, at the limit $q \to 1$ we have $z^*z=zz^*$. This allows one
to treat ${\rm Pol}({\Bbb C)}_q$ as a polynomial algebra on the quantum
plane.

 It is worthwhile to note that the passage to the generators
$a^+=(1-q^2)^{-1/2}z^*$, $a=(1-q^2)^{-1/2}z$ realizes an isomorphism between
the $*$-algebra ${\rm Pol}({\Bbb C})_q$ and the q-analogue of the Weil
algebra considered in \cite{CK}: $a^+a-q^2aa^+=1$.

 Impose the notation $y=1-zz^*$. It is straightforward that
$$z^*y=q^2yz^*,\quad zy=q^{-2}yz.\eqno(1.1)$$
It is also easy to show that any element $f \in{\rm Pol}({\Bbb C)}_q$ admits
a unique decomposition
$$f=\sum_{m>0}z^m
\psi_m(y)+\psi_0(y)+\sum_{m>0}\psi_{-m}(y)z^{*m}.\eqno(1.2)$$
The passage to the decomposition (1.2) is similar to the passage from
Cartesian coordinates on the plane ${\Bbb C}\simeq{\Bbb R}^2$ to polar
coordinates.

 We follow \cite{KL, NN} in completing the vector space ${\rm Pol}({\Bbb
C)}_q$ to obtain the function space in the quantum disc.

 Consider the ${\rm Pol}({\Bbb C)}_q$-module $H$ determined by its generator
$v_0$ and the relation $z^*v_0=0$. Let $T:{\rm Pol}({\Bbb C)}_q \to{\rm
End}_{\Bbb C}(H)$ be the representation of ${\rm Pol}({\Bbb C)}_q$ in $H$,
and $l_{m,n}(f)$, $m,n \in {\Bbb Z}_+$, stand for the matrix elements of the
operator $T(f)$ in the basis $\{z^mv_0\}_{m \in {\Bbb Z}_+}$.

 Impose three topologies in the vector space ${\rm Pol}({\Bbb C)}_q$ and
prove their equivalence.

 Let ${\cal T}_1$ be the weakest among the topologies in which all the
linear functionals $l_{j,k}^\prime:f \mapsto a_{jk}(f)$, $j,k \in{\Bbb
Z}_+$, are continuous,\\
${\cal T}_2$ the weakest among the topologies in which all the linear
functionals $l_{m,n}^{\prime \prime}:f \mapsto \psi_m(q^{2n})$, $m \in{\Bbb
Z}$, $n \in{\Bbb Z}_+$, are continuous,\\
${\cal T}_3$ the weakest among the topologies in which all the linear
functionals $l_{m,n}$, $m,n \in{\Bbb Z}_+$, are continuous.

\medskip

\begin{proposition}\hspace{-.5em}. The topologies ${\cal T}_1, {\cal T}_2,
{\cal T}_3$ are equivalent.
\end{proposition}

\smallskip

 {\bf Proof.} Remind the standard notation:
$$(t;q)_m=\prod_{j=0}^{m-1}(1-tq^j);\quad
(t;q)_\infty=\prod_{j=0}^\infty(1-tq^j).$$
It follows from the definitions that with $m \ge k$
$$z^{*k}z^mv_0=(q^{2m};q^{-2})_k \cdot z^{m-k}v_0.$$
Hence
$$l_{mn}(f)=\sum_{j=0}^{\min(m,n)}(q^{2m};q^{-2})_{m-j} \cdot
l_{n-j,m-j}^\prime;\quad m,n \in{\Bbb Z}_+.$$
Thus the topology ${\cal T}_1$ is stronger than ${\cal T}_3$. In fact, these
topologies are equivalent since the linear span of the functionals
$\{l_{m-j,n-j}\}_{j=0}^{\min(m,n)}$ coincides with that of
$\{l_{m-j,n-j}^\prime \}_{j=0}^{\min(m,n)}$. The equivalence of ${\cal T}_2$
and ${\cal T}_3$ follows from the relations $l_{jk}=(q^{2j};q^{-2})_{\max(0,j-k)}l_{k-j,\min(j,k)}^{\prime \prime}$.
\hfill $\Box$

\medskip

 We denote the completion of the Hausdorff topological vector space ${\rm
Pol}({\Bbb C)}_q$ by $D(U)_q^\prime$ and call it the space of distributions
in the quantum disc. $D(U)_q^\prime$ may be identified with the space of
formal series of the form (1.2) whose coefficients $\psi_j(y)$ are defined
in $q^{2{\Bbb Z}_+}$. The linear functionals $l_{mn}$ are extendable by
continuity onto the topological vector space $D(U)_q^\prime$. Associate to
each distribution $f \in D(U)_q^\prime$ the infinite matrix
$T(f)=(l_{mn}(f))_{m,n \in{\Bbb Z}_+}$.

 A distribution $f \in D(U)_q^\prime$ is said to be finite if
$\#\{(j,k)|\,\psi_j(q^{2k})\ne 0\}<\infty$. Evidently, a distribution $f$ is
finite iff the matrix $T(f)$ has only finitely many non-zero entries. The
vector space of finite functions in the quantum disc is denoted by
$D(U)_q$. There exists a non-degenerate pairing
$$D(U)_q^\prime \times D(U)_q \to {\Bbb C};\quad f_1 \times f_2 \mapsto
{\rm tr}\,T(f_1)T(f_2).$$

 The extension by continuity procedure allows one to equip $D(U)_q$ with a
structure of  $^*$-algebra, and $D(U)_q^\prime$ with a structure of
$D(U)_q$-bimodule.

 Consider the algebra ${\rm Pol}({\Bbb C})_q^{op}$ derived from ${\rm
Pol}({\Bbb C})_q$ via a replacement of the multiplication law by the
opposite one. The elements $z \otimes 1,\,z^*\otimes 1,\,1 \otimes z,\,1
\otimes z^*$ of ${\rm Pol}({\Bbb C})_q^{op}\otimes{\rm Pol}({\Bbb C})_q$ are
denoted respectively by $z,z^*,\zeta,\zeta^*$. To avoid confusion in the
notation, we use braces to denote the multiplication in ${\rm Pol}({\Bbb
C})_q^{op}\otimes{\rm Pol}({\Bbb C})_q$, e.g. $\{zz^*\}=q^2\{z^*z\}+1-q^2$.
The module $H^{op}$ over ${\rm Pol}({\Bbb C})_q^{op}$ is defined by its
generator $v_0^{op}$ and the relation $zv_0^{op}=0$. Apply the above
argument to ${\rm Pol}({\Bbb C})_q^{op}\otimes{\rm Pol}({\Bbb C})_q$-module
$H^{op}\otimes H$ in order to introduce the algebra $D(U \times U)_q$ of
finite functions in the Cartesian product of quantum discs, together with
$D(U \times U)_q$-bimodule $D(U \times U)_q^\prime$.

 The reason for replacement the multiplication law in ${\rm Pol}({\Bbb
C})_q$ by the opposite one will become clear later (see \cite{SSV2, V1, V2}).

 The linear functional
$$\mu(f)=(1-q^2)\sum_{m \in{\Bbb Z}_+}\psi_0(q^{2m})q^{2m},\quad f \in D(U)_q,$$
is called the (normalized) Lebesgue integral in the quantum disc, since
under the formal passage to the limit as $q \to 1$ one has $\mu(f)\to \frac{1}{\pi}
\begin{array}[t]{c}\int \int \\ \scriptstyle U \end{array}
fd{\rm Im}\,z \cdot d{\rm Re}\,z$.

 Let $K \in D(U \times U)_q^\prime$; the integral operator $f \mapsto {\rm
id}\otimes \mu(K(1 \otimes f))$ with kernel $K$ maps $D(U)_q$ into
$D(U)_q^\prime$. We are interested in solving an inverse problem which is in
finding out the explicit formulae for kernels $K \in D(U \times U)_q^\prime$
of well known linear operators. In this field, an analogue of the Bergman
kernel for the quantum disc was obtained in the work of S. Klimek and A.
Lesniewski \cite{KL}:
$$K_q(z,\zeta)=(1-z \zeta^*)^{-1}(1-q^2z \zeta^*)^{-1}.$$

 Finally, equip $H$ with the structure of a pre-Hilbert space by setting
$$(z^jv_0,z^mv_0)=\delta_{jm}(q^2;q^2)_m;\quad j,m \in{\Bbb Z}_+.$$
It is easy to show that $T(z^*)=T(z)^*$, $I-T(z)T(z^*)\ge 0$. Thus we get a
$*$-representation of ${\rm Pol}({\Bbb C})_q$ in the completion
$\overline{H}$ of $H$ (see \cite{KL}).

\bigskip

\section{Differential forms and $\overline{\partial}$-problem}

 Let $\Omega_q({\Bbb C})$ stand for the involutive algebra determined by its
generators $z,dz$ and the relations
$$1-z^*z=q^2(1-zz^*),\quad dz \cdot z^*=q^{-2}z^*\cdot dz,\quad dz \cdot
z=q^2z \cdot dz,$$
$$dz \cdot dz^*=-q^{-2}dz^*dz,\quad dz \cdot dz=0.$$
(Also, an application of the involution $^*$ to the above yields
$$dz^* \cdot dz^*=0,\quad dz^*\cdot z=q^2z \cdot dz^*,\quad dz^* \cdot
z^*=q^{-2}z^*\cdot dz^*.)\eqno(2.1)$$
Equip $\Omega_q({\Bbb C})$ with the grading as follows:
$${\rm deg}\,z={\rm deg}\,z^*=0,\quad{\rm deg}\,dz={\rm deg}\,dz^*=1.$$
There exists a unique linear map $d:\Omega_q({\Bbb C})\to \Omega_q({\Bbb
C})$ such that
$$d:z \mapsto dz,\quad d:z^*\mapsto dz^*,\quad d:dz \mapsto 0,\quad
d:dz^*\mapsto 0,$$
and
$$d(\omega^\prime \cdot \omega^{\prime \prime})=d \omega^\prime \cdot
\omega^{\prime \prime}+(-1)^{{\rm deg}\,\omega^\prime}\cdot \omega^\prime d
\omega^{\prime \prime};\quad \omega^\prime,\omega^{\prime \prime}\in
\Omega_q({\Bbb C}).$$
Evidently, $d^2=0$, and $(d \omega)^*=d \omega^*$ for all $\omega \in
\Omega_q({\Bbb C})$.

 Turn to a construction of operators $\partial,\overline{\partial}$. For
that, we need a bigrading in $\Omega_q({\Bbb C})$:
$${\rm deg}\,z={\rm deg}\,z^*=(0,0);\quad {\rm deg}\,(dz)=(1,0);\quad{\rm
deg}\,(dz^*)=(0,1).$$
Now $d$ has a degree 1 and admits a unique decomposition into a sum
$d=\partial+\overline{\partial}$ of operators $\partial$,
$\overline{\partial}$ with bidegrees respectively (1,0) and (0,1). A standard
argument allows one to deduce from $d^2=0$ that
$\partial^2=\overline{\partial}^2=\partial
\overline{\partial}+\overline{\partial}\partial=0$. It is also easy to show
that $(\partial \omega)^*=\overline{\partial}\omega^*$ for all $\omega \in
\Omega_q({\Bbb C})$.

 Each element $\omega \in \Omega_q({\Bbb C})$ is uniquely decomposable into
a sum
$$\omega=f_{00}+dz \,f_{10}+f_{01}dz^*+dz \,f_{11}dz^*,\quad f_{ij}\in
{\rm Pol}({\Bbb C})_q,\quad i,j=0,1.$$
Equip $\Omega_q({\Bbb C})$ with a topology corresponding to this
decomposition:
$$\Omega_q({\Bbb C})\simeq{\rm Pol}({\Bbb C})_q \oplus{\rm Pol}({\Bbb C})_q
\oplus{\rm Pol}({\Bbb C})_q \oplus{\rm Pol}({\Bbb C})_q.$$
Pass as above via a completion procedure from ${\rm Pol}({\Bbb C})_q$ to the
space of distributions $D(U)_q^\prime$ and then to the space of finite
functions to obtain the bigraded algebra $\Omega(U)_q$. The operators
$d,\partial,\overline{\partial}$ are transferred by continuity from
$\Omega_q({\Bbb C})$ onto the algebra $\Omega(U)_q$ of differential forms
with finite coefficients in the quantum disc.

 The subsequent constructions involve essentially q-analogues of type (0,*)
differential forms with coefficients in sections of holomorphic bundles. The
latter carry a structure of bimodules over algebras of type (0,*)
differential forms as above. Remind the notion of a differentiation for such
a bimodule.

 Let $\Omega$ be a ${\Bbb Z}_+$-graded algebra and $M$ a ${\Bbb Z}_+$-graded
$\Omega$-bimodule. A degree 1 operator is said to be a differentiation if
for all $m \in M,\,\omega \in \Omega$ one has $\overline{\partial}(m
\omega)=(\overline{\partial}m)\omega+(-1)^{{\rm deg}\,m}m \cdot
\overline{\partial}\omega$, $\overline{\partial}(\omega
m)=(\overline{\partial}\omega)\cdot m+(-1)^{{\rm deg}\,\omega}\omega \cdot
\overline{\partial}m$.

 Let $\lambda \in {\Bbb C}$. Consider the graded bimodule over $\Omega({\Bbb
C})_q^{(0,*)}=\Omega({\Bbb C})_q^{(0,0)}+\Omega({\Bbb C})_q^{(0,1)}$
determined by its generator $v_ \lambda$ with ${\rm deg}(v_ \lambda)=0$ and
the relations
$$z \cdot v_ \lambda=q^{-\lambda}v_ \lambda \cdot z,\quad z^*\cdot v_
\lambda=q^{\lambda} v_ \lambda \cdot z^*,\quad dz^*\cdot v_ \lambda=q^{\lambda}
v_ \lambda \cdot dz^*.$$
We denote this graded bimodule by $\Omega({\Bbb C})_{\lambda,q}^{(0,*)}$. It
possesses a unique degree 1 differentiation $\overline{\partial}$ such that
$\overline{\partial}v_ \lambda=0$. Pass (via an extension by continuity)
from polynomial coefficients to finite ones to obtain the graded bimodule
$\Omega(U)_{\lambda,q}^{(0,*)}$ over $\Omega(U)_q^{(0,*)}$, together with its
differentiation $\overline{\partial}$.

 We restrict ourselves to the case $\lambda \in{\Bbb R}$ and equip the spaces
$\Omega(U)_{\lambda,q}^{(0,0)},\,\Omega(U)_{\lambda,q}^{(0,1)}$ with the scalar products
$$(f_1 \cdot v_ \lambda,f_2 \cdot v_
\lambda)=\int_{U_q}f_2^*f_1(1-zz^*)^{\lambda-2}d \mu,\eqno(2.4)$$
$$(f_1v_ \lambda dz^*,f_2v_
\lambda dz^*)=\int_{U_q}f_2^*f_1(1-zz^*)^{\lambda}d \mu.\eqno(2.5)$$

 The completions of the pre-Hilbert spaces $\Omega(U)_{\lambda,q}^{(0,0)}$ and
$\Omega(U)_{\lambda,q}^{(0,1)}$ can be used in the formulation of
$\overline{\partial}$-problem. Specifically, we mean finding a solution of
the equation $\overline{\partial}u=f$ in the orthogonal complement to the
kernel of $\overline{\partial}$. In the classical case ($q=1$) such a
formulation is standard \cite{BFG}, and the solution is very well known. If
$\lambda$ stand for a real number and $\lambda>1$ than one has $$u(z)={1
\over 2 \pi i}\int_U{1 \over z-\zeta}\left({1-|\zeta|^2 \over
1-\overline{\zeta}z}\right)^{\lambda-1}f(\zeta)d \overline{\zeta}\wedge d
\zeta.\eqno(2.6)$$

 (2.6) implies the "Cauchy-Green formula": $$u(z)={{\lambda-1} \over 2 \pi
i}\int_U{(1-|\zeta|^2)^{\lambda-2}
\over(1-\overline{\zeta}z)^{\lambda}}u(\zeta)d \overline{\zeta}\wedge d
\zeta+{1 \over 2 \pi i}\int_U{(1-|\zeta|^2)^{\lambda-1}
\over(z-\zeta)(1-\overline{\zeta}z)^{\lambda-1}}{\partial u \over \partial
\overline{\zeta}}d \overline{\zeta}\wedge d \zeta.\eqno(2.7)$$. Our purpose
is to obtain the q-analogues of (2.6), (2.7) for $\lambda=2$.

 The standard way of solving the $\overline{\partial}$-problem is to solve
first the Poisson equation $\Box \omega=f$ with
$\Box=-\overline{\partial}^*\overline{\partial}$. In the case $\lambda=2$
the kernel in (2.6) is derived by a differentiation in $z$ of the Green
function $$G(z,\zeta)={1 \over \pi}{\rm ln}(|z-\zeta|^2)-{1 \over \pi}{\rm
ln}(|1-z \overline{\zeta}|^2).\eqno(2.8)$$ In its turn, (2.8) can be
obtained by the d'Alembert method: the first term is contributed by a real
source, and the second one is coming from an imaginary source.

 Note that the differential calculus for the quantum disc we use here is
well known (see, for example, \cite{Zu}). Its generalization onto the case
of an arbitrary bounded symmetric domain was obtained in \cite{SV} via an
application of a quantum analogue of the Harish-Chandra embedding \cite{H1}.

\bigskip

\section{Green function for Poisson equation}

 With $q=1$ the measure $d \nu=(1-|z|^2)^{-2}d \mu$ is invariant with
respect to the M\"obius transformations. In the case $q \in(0,1)$ impose an
"invariant integral" $\nu:D(U)_q \to{\Bbb C}$, $f \mapsto \int
\limits_{U_q}fd \nu$ by setting
$$\int \limits_{U_q}fd \nu \stackrel{def}{=}\int \limits_{U_q}f
\cdot(1-zz^*)^{-2}d \mu.$$
 The Hilbert spaces $L^2(d \nu)_q$, $L^2(d \mu)_q$ are defined as
completions of the vector spaces $D(U)_q=\Omega(U)_q^{(0,0)}$,
$\Omega(U)_q^{(0,1)}$ with respect to the norms
$$\Vert f \Vert=\left(\int \limits_{U_q}f^*fd \nu \right)^{1/2}, \quad \Vert
fdz^* \Vert=\left(\int \limits_{U_q}f^*fd \mu \right)^{1/2}.$$

  Proofs of the following statements are to be found in
\cite[Proposition 5.7, Corollary 5.8]{SSV2}, \cite [Corollary 4.2]{SSV4}.

\medskip

\begin{lemma}\hspace{-.5em}. There exist $0<c_1 \le c_2$ such that
$$c_1 \le \overline{\partial}^*\overline{\partial}\le c_2.\eqno(3.1)$$
\end{lemma}

\medskip

\begin{proposition}\hspace{-.5em}. The exact estimates for
$\overline{\partial}^*\overline{\partial}$ are of the form
$${1 \over (1+q)^2}\le \overline{\partial}^*\overline{\partial}\le
{1 \over (1-q)^2}.\eqno(3.2)$$
\end{proposition}

\medskip

 The inequalities (3.1) allow one to extend by continuity the operators
$\overline{\partial}$, $\Box=-\overline{\partial}^*\overline{\partial}$ from
the dense subspace of finite functions $D(U)_q$ onto the entire $L^2(d
\nu)_q$. They also imply that for any $f \in L^2(d \nu)_q$ there exists a
unique solution $u$ of Poisson equation $\Box u=f$. Now it follows from
(3.2) that $\Vert u \Vert \le (1+q)^2 \Vert f \Vert$.

 One obtains easily from definitions

\medskip

\begin{lemma}\hspace{-.5em}. The series
$$\sum_{i,j \in{\Bbb Z}_+}z^{*i} \zeta^{i}\psi_{ij}(y,\eta)z^{}\zeta^{*j},\quad
y=1-z^*z,\quad \eta=1-\zeta\zeta^*$$
converges in $D(U \times U)_q^\prime$ for any family
$\{\psi_{ij}(y,\eta)\}_{i,j \in{\Bbb Z}_+}$ of functions defined on $q^{2{\Bbb
Z}_+}\times q^{2{\Bbb Z}_+}$.
\end{lemma}

\medskip

\begin{corollary}\hspace{-.5em}. For all $m \ge 0$ there exists a well
defined generalized kernel
$$G_m=\left \{\left((1-\zeta \zeta^*)(1-z^*\zeta)^{-1}\right)^m
\left((1-z^*z)(1-z \zeta^*)^{-1}\right)^m \right \}.\eqno(3.3)$$
\end{corollary}

\medskip

 To state the principal result of the section, we need an expansion of the
Green function (2.8):
$${\rm ln}{|z-\zeta|^2 \over|1-z \overline{\zeta}|^2}={\rm
ln}\left(1-\frac{(1-|z|^2)(1-|\zeta|^2)}{|1-z
\overline{\zeta}|^2}\right)=-\sum_{m=1}^\infty{1 \over
m}\left(\frac{(1-|z|^2)(1-|\zeta|^2)}{|1-z \overline{\zeta}|^2}\right)^m.$$

 Evidently, a formal passage to a limit yields
$$\lim_{q \to 1}\,G_m=\left(\frac{(1-|z|^2)(1-|\zeta|^2)}{|1-z
\overline{\zeta}|^2}\right)^m.$$

 A proof of the following result one can find in \cite[Theorem 1.2]{SSV4}.

\medskip

\begin{theorem}\hspace{-.5em}. The continuous operator $\Box^{-1}$ in $L^2(d
\nu)_q$ coincides on the dense linear subspace $D(U)_q \subset L^2(d \nu)_q$
with the integral operator whose kernel is $G=- \displaystyle \sum
\limits_{m=1}^\infty{\textstyle q^{-2}-1 \over \textstyle q^{-2m}-1}G_m$:
$$\Box^{-1}f=\int_{U_q}G(z,\zeta)f(\zeta)d \nu.$$
Here $G_m \in D(U \times U)_q^\prime$ is given by (3.3).
\end{theorem}

\medskip

 Note in conclusion that the operators $\partial$, $\overline{\partial}$,
$\Box$ admit an extension by continuity onto the space $D(U)_q^\prime$ of
distributions in the quantum disc.

\bigskip

\section{Cauchy-Green formula}

 One can use the differentials $\partial:\Omega_q^{(0,0)}\to
\Omega_q^{(1,0)}$, $\overline{\partial}:\Omega_q^{(0,0)}\to
\Omega_q^{(0,1)}$ to define the partial derivatives ${\textstyle
\partial^{(l)}\over \textstyle \partial z}$, ${\textstyle
\partial^{(r)}\over \textstyle \partial z}$, ${\textstyle
\partial^{(l)}\over \textstyle \partial z^*}$, ${\textstyle
\partial^{(r)}\over \textstyle \partial z^*}$. Specifically, we set up
$\partial f=dz \cdot{\textstyle \partial^{(l)}f \over \textstyle \partial
z}={\textstyle \partial^{(r)}f \over \textstyle \partial z}dz$,
$\overline{\partial}f=dz^* \cdot{\textstyle \partial^{(l)}f \over \textstyle
\partial z^*}={\textstyle \partial^{(r)}f \over \textstyle \partial
z^*}dz^*$. It is easy to show that these operators admit extensions by
continuity from $D(U)_q$ onto $D(U)_q^\prime$.

 Let $f \in D(U)_q$. Define the integral of the (1,1)-form $dz \cdot f \cdot
dz^*$ over the quantum disc by $\int \limits_{U_q}dz \cdot f \cdot dz^*=-2i \pi
\int \limits_{U_q}fd \mu$.

 A proof of the following proposition can be found in \cite[Theorem
2.1]{SSV4}.

\medskip

\begin{proposition}\hspace{-.5em}. Let $f \in D(U)_q$. Then
\begin{enumerate}

\item There exists a unique solution $u \in L^2(d \mu)_q$ of the
$\overline{\partial}$-problem $\overline{\partial}u=f$, which is orthogonal
to the kernel of $\overline{\partial}$.

\item $u= {\textstyle 1 \over \textstyle 2 \pi i}\displaystyle \int \limits_{U_q} d \zeta {\textstyle \partial^{(l)}\over
\textstyle \partial z}G(z,\zeta) fd \zeta^*$, with $G
\in D(U \times U)_q^\prime$ being the Green function of the Poisson
equation.

\item $f=-{\textstyle 1 \over \textstyle 2 \pi i}\displaystyle \int
\limits_{U_q}(1-z \zeta^*)^{-1}(1-q^{-2}z \zeta^*)^{-1}d \zeta f(\zeta)d
\zeta^*- {\textstyle 1 \over \textstyle 2 \pi i}\displaystyle \int
\limits_{U_q} d \zeta{\textstyle \partial^{(l)}\over \textstyle \partial
z}G(z,\zeta) \cdot{\textstyle \partial^{(r)}f \over \textstyle \partial
\zeta^*}d \zeta^*$.

\end{enumerate}
\end{proposition}

\bigskip

\section{Eigenfunctions of the operator $\Box$}

 Let ${\Bbb C}[\partial U]_q$ stand for the algebra of finite sums of the
form
$$\sum_{m \in{\Bbb Z}}a_me^{im \theta},\quad \theta \in {\Bbb R}/2 \pi{\Bbb
Z}\eqno(5.1)$$
with complex coefficients. The ${\Bbb C}[\partial U]_q$-module of formal
series of the form (5.1) is denoted by ${\Bbb C}[[\partial U]]_q$. We also
denote the algebra of finite sums like (5.1) with coefficients from $D(U)_q$
by $D(U \times \partial U)_q$, and the module of formal series (5.1) with
coefficients from $D(U)_q^\prime$ by $D(U \times \partial U)_q^\prime$. This
vector space will be equipped by the topology of coefficientwise convergence.

 The use of the index $q$ in the notation for the above vector spaces is
justified by the fact that, as one can show, they are in fact modules over
the quantum universal enveloping algebra.

 Recall the notations \cite{GR}:
$$(a;q^2)_\infty=\prod_{j \in{\Bbb Z}_+}(1-aq^{2j}),\quad(a;q^2)_\gamma=
\frac{(a;q^2)_\infty}{(aq^{2 \gamma};q^2)_\infty},\quad \gamma \in{\Bbb C}.$$

 With $q=1$, the integral
$$u(z)=\int \limits_{\partial U}\left(\frac{1-|z|^2}{(1-z \overline
\zeta)(1-\overline z \zeta)}\right)^{l+1}f(\zeta)d \nu,\quad d \nu={d \theta
\over 2 \pi},$$
represents an eigenfunction of $\Box$ (see \cite{H2}):
$$\Box u=\lambda(l)u,\quad \lambda(l)=\left(l+{1 \over 2}\right)^2-{1 \over
4}.$$

 With $q \in(0,1)$, the power $P^\gamma$ of the Poisson kernel
$P={\textstyle 1-|z|^2 \over \textstyle|1-z \overline \zeta|^2}$ is replaced
by the element $P_ \gamma \in D(U \times \partial U)_q$:
$$P_ \gamma=(1-zz^*)^\gamma(z
\zeta^*;q^2)_{-\gamma}\cdot(q^2z^*\zeta;q^2)_{-\gamma}.\eqno(5.2)$$

 Here $(z \zeta^*;q^2)_{-\gamma}$, $(q^2z^*\zeta;q^2)_{-\gamma}$ are the
$q$-analogues of the powers $(1-z \overline \zeta)^{-\gamma}$, $(1-\overline
z \zeta)^{-\gamma}$, and the $q$-binomial theorem (see \cite{GR}) implies
$$(z \zeta^*;q^2)_{-\gamma}=\sum_{n \in{\Bbb Z}_+}{(q^{2 \gamma};q^2)_n
\over(q^2;q^2)_n}(q^{-2 \gamma}z \zeta^*)^n,$$
$$(q^2z^*\zeta;q^2)_{-\gamma}=\sum_{n \in{\Bbb Z}_+}{(q^{2 \gamma};q^2)_n
\over(q^2;q^2)_n}(q^{2-2 \gamma}z^*\zeta)^n.$$

 The following proposition is proved in \cite[Theorem 3.1]{SSV4}.

\medskip

\begin{proposition}\hspace{-.5em}. For all $f \in{\Bbb C}[\partial U]_q$ the
element
$$u=\int \limits_{\partial U}P_{l+1}(z,e^{i \theta})f(e^{i \theta}){d \theta
\over 2 \pi}\eqno(5.3)$$
of $D(U)_q^\prime$ is an eigenvector of $\Box$:
$$\Box u=\lambda(l)u,\quad
\lambda(l)=-\frac{(1-q^{-2l})(1-q^{2l+2})}{(1-q^2)^2}.$$
\end{proposition}

\medskip

 We need the following standard notation (\cite{GR}):
$$_r \Phi_s \left[{a_1,a_2,\ldots,a_r;q;z \atop b_1,\ldots,b_s}\right]=$$
$$=\sum_{n \in{\Bbb Z}_+} \frac{(a_1;q)_n \cdot(a_2;q)_n \cdot \ldots
\cdot(a_r;q)_n} {(b_1;q)_n \cdot(b_2;q)_n \cdot \ldots
\cdot(b_s;q)_n(q;q)_n}\left((-1)^n \cdot q^{n(n-1)\over
2}\right)^{1+s-r}\cdot z^n.$$

\medskip

\begin{corollary} (cf. \cite{V1}). The series
$$\varphi_l(y)=\,_3 \Phi_2 \left[{y^{-1},q^{-2l},q^{2(l+1)};q^2;q^2 \atop
q^2;0}\right]$$
converges in $D(U)_q^\prime$, and its sum is an eigenfunction of $\Box$:\
$\Box \varphi_l=\lambda(l)\varphi_l$.
\end{corollary}

\smallskip

 {\bf Proof}. The convergence of the series is due to the fact that it
breaks for each $y \in q^{2{\Bbb Z}_+}$. So it suffices to establish the
relation
$$\varphi_l(y)=\int \limits_{\partial U}P_{l+1}(z,\zeta)d \nu.$$

 It follows from the definitions that the above integral equals to
$$\sum_{n \in{\Bbb Z}_+}
\frac{(q^{2l+2};q^2)_n^2}{(q^2;q^2)_n^2}q^{-2(2l+1)n}y^{l+1}z^nz^
{*n}=y^{l+1}\,
_3 \Phi_1 \left[{q^{2+2l},q^{2+2l},y^{-1};q^2;q^{-2(2l+1)}y \atop
q^2}\right].$$

 Now it remains to apply the identity (see \cite{GR}):
$$b^n \,_3 \Phi_1 \left[{q^{-n},b,{\textstyle q \over
\textstyle z};q,{\textstyle z \over \textstyle c}\atop {\textstyle
bq^{1-n}\over \textstyle c}}\right]=\,_3 \Phi_2 \left[{q^{-n},b,{\textstyle
bzq^{-n}\over \textstyle c};q,q \atop {\textstyle bq^{1-n}\over \textstyle
c},0}\right],$$
with $q$ being replaced by $q^2$, $y$ by $q^{2n}$, $b$ by $q^{2l+2}$, $z$ by
$q^{-2l}$, and $c$ by $q^{2+2l-2n}$. \hfill $\Box$.

 Note that $\varphi_l(y)$ is a $q$-analogue of a spherical function on a
hyperbolic plane (see \cite{H2}).

 For each $l \in{\Bbb C}$ a linear operator has been constructed from ${\Bbb
C}[\partial U]_q$ into the eigenspace of $\Box$, associated to the
eigenvalue $\lambda(l)$. Now we try to invert this linear operator.

 For that, we need a $q$-analogue of the operator $b_r:f(z)\mapsto f(re^{i
\theta})$ which restricts the function in the disc onto the circle $|z|=r$
of radius $r \in(0,1)$. Let $r>0$, $1-r^2 \in q^{2{\Bbb Z}_+}$. Define a
linear operator $b_r:D(U)_q^\prime \to{\Bbb C}[[\partial U]]_q$ by
$$b_r:\sum_{j>0}z^j \cdot \psi_j(y)+\psi_0(y)+\sum_{j>0}\psi_{-j}(y) \cdot
z^{*j}\mapsto$$
$$\sum_{j>0}(re^{i \theta})^j \cdot
\psi_j(1-r^2)+\psi_0(1-r^2)+\sum_{j>0}\psi_{-j}(q^{-2j}\cdot (1-r^2))(re^{-i
\theta})^j.$$
(It is implicit that the functions $\psi_j(y)$, $j \in{\Bbb Z}$, vanish at $y
\notin q^{2{\Bbb Z}_+})$.

 Recall the definition of the $q$-gamma-function (\cite{GR}):
$$\Gamma_q(x)={(q;q)_\infty \over (q^x;q)_\infty}(1-q)^{1-x}.$$

 One may assume without loss of generality that
$$0 \le{\rm Im}\,l<{\pi \over 2ln(q^{-1})},\quad {\rm Re}\,l \ge-{1 \over
2}.$$

\medskip

\begin{proposition}\hspace{-.5em}. Let ${\rm Re}\,l>-{\textstyle 1 \over
\textstyle 2}$, and $u \in D(U)_q^\prime$ is an eigenfunction of $\Box$
given by (5.3). Then in ${\Bbb C}[\partial U]_q$ one has
$$f=\frac{\Gamma_{q^2}^2(l+1)}{\Gamma_{q^2}(2l+1)}\lim_{1-r^2 \in q^{2{\Bbb
Z}_+},\,r \to 1}(1-r^2)^lb_ru.$$
\end{proposition}

\medskip

 The proof of this proposition is based on the following result which was
communicated to the authors by L. I. Korogodsky:

\medskip

\begin{lemma}\hspace{-.5em}.
$$\lim \limits_{x \in q^{-2{\Bbb Z}_+},\,x \to \infty}\varphi_l \left({1
\over x}\right)\left/
\left(\frac{\Gamma_{q^2}(2l+1)}{\Gamma_{q^2}^2(l+1)}x^l \right)=1
\right.\leqno 1).$$
if ${\rm
Re}\,l>-{\textstyle 1 \over \textstyle 2}$.

$$\lim \limits_{x \in q^{-2{\Bbb Z}_+},\,x \to \infty}\varphi_l \left({1
\over x}\right)\left/
\left(\frac{\Gamma_{q^2}(-2l-1)}{\Gamma_{q^2}^2(-l)}x^{-l-1}\right)=1
\right.\leqno 2).$$
if ${\rm Re}\,l<-{\textstyle
1 \over \textstyle 2}$.
\end{lemma}

\smallskip

 {\bf Proof}. It follows from the relation $\varphi_l(y)=\varphi_{-1-l}(y)$
that one may restrict oneself to the case ${\rm Re}\,l>-{\textstyle 1 \over
\textstyle 2}$. An application of the identity (\cite{GR})
$$_2 \Phi_1 \left[{q^{-n},b;q;z \atop c}\right]=\frac{\left({\textstyle c
\over \textstyle b};q \right)_n}{(c;q)_n}\;_3 \Phi_2
\left[{q^{-n},b,{\textstyle bzq^{-n}\over \textstyle c};q;q \atop {\textstyle
bq^{1-n}\over \textstyle c},0}\right],$$
with $q$, $b$, $c$, $z$ being replaced respectively by $q^2$, $q^{-2l}$,
$q^{-2l-2n}$, $q^{2l+2}$, yields
$$\varphi_l(q^{2n})=\frac{(q^{-2l-2n};q^2)_n}{(q^{-2n};q^2)_n}\cdot \,_2
\Phi_1 \left[{q^{-2n};q^{-2l};q^2;q^{2l+2} \atop q^{-2l-2n}}\right]\sim$$
$$\sim q^{-2nl}\frac{(q^{2(l+1)};q^2)_\infty}{(q^2;q^2)_\infty}\cdot \,_1
\Phi_0[q^{-2l};q^2;q^{2(l+1)}]=$$
$$=q^{-2nl}\frac{(q^{2(l+1)};q^2)_\infty}{(q^2;q^2)_\infty}\cdot
\frac{(q^{2(l+1)};q^2)_\infty}{(q^{2(2l+1)};q^2)_\infty}.$$
Now it remains to refer to the definition of the $q$-gamma-function. \hfill
$\Box$

\medskip

 In the special case $f=1$ proposition 5.3 follows from lemma 5.4. The
general case reduces to the above special case via an application of a
quantum symmetry argument, which will be described in \cite[Theorem
3.7]{SSV4}.

\bigskip

\section{Fourier transformation}

 It follows from the definitions that the integral operator with kernel
$K=\sum \limits_ik_i^{\prime \prime}\otimes k_i^\prime$ is conjugate to the
integral operator with the kernel $K^t=\sum \limits_ik_i^{\prime*}\otimes
k_i^{\prime \prime*}$. Note that the conjugate to the unitary is an inverse
operator.

 Recall \cite{H2} the heuristic argument that leads to the Fourier
transformation. Proposition 5.1 allows one to obtain eigenfunctions of
$\Box$. It is natural to expect that "any" function $u$ admits a
decomposition in eigenfunctions of $\Box$, and that the associated Fourier
operator is unitary.

 Impose the notations: $h=-2 \ln q$,
$$P_{l+1}^t=(q^2z^*\zeta;q^2)_{-l-1}(z \zeta^*;q^2)_{-l-1}(1-\zeta
\zeta^*)^{1+l},$$
$$c(l)=\Gamma_{q^2}(2l+1)/(\Gamma_{q^2}(l+1))^2.$$

 It is shown in \cite[section 5]{SSV4} that, just as in the standard
representation theory (see \cite{H2}), one has

\medskip

\begin{proposition}\hspace{-.5em}. Consider the Borel measure $d \sigma$ on
$[0,{\textstyle  \pi \over \textstyle h}]$, given by
$$d \sigma(\rho)={1 \over 2 \pi}\cdot {h \cdot e^h \over e^h-1}c(-{1 \over
2}+i \rho)^{-1}\cdot c(-{1 \over 2}-i \rho)^{-1}d \rho.$$ The integral
operators $$u(z)\mapsto \int \limits_{U_q}P_{{1 \over 2}-i
\rho}^t(z,\zeta)u(\zeta)d \nu,$$ $$f(e^{i \theta},\rho)\mapsto \int
\limits_0^{ \pi/h}\int \limits_0^{2 \pi}P_{{1 \over 2}+i \rho}(z,e^{i
\theta})f(e^{i \theta},\rho){d \theta \over 2 \pi}d \sigma(\rho)$$ are
extendable by continuity from the dense linear subspaces $${\rm Pol}({\Bbb
C})_q \subset L^2(d \nu)_q,\quad C^\infty[0,{\textstyle  \pi \over
\textstyle h}]\otimes{\Bbb C}[\partial U]_q \subset L^2({d \theta \over 2
\pi})\otimes L^2(d \sigma)$$ up to mutually inverse unitaries $F$, $F^{-1}$.
\end{proposition}

\medskip

 {\bf Remark 6.2}. The function $c(l)$, the measure $d \sigma(\rho)$ and the
operator $F$ are the quantum analogues for c-function of Harish-Chandra,
Plancherel measure and Fourier transformation respectively (see \cite{H2}).

\bigskip

\section{Berezin deformation of the quantum disc}

 We are going to use in the sequel bilinear operators $L:D(U)_q \times D(U)_q
\to D(U)_q$ of the form
$$L:\,f_1 \times f_2 \to
\sum_{ijkm=0}^{N(L)}a_{ijkm}\left(\left({\partial^{(r)}\over \partial
z^*}\right)^if_1 \right)z^{*j}z^k \left(\left({\partial^{(l)}\over \partial
z}\right)^mf_2 \right),\eqno(7.1)$$
with $a_{ijkm}\in{\Bbb C}$. Such operators will be called $q$-bidifferential.

 Our principal purpose is to construct the formal deformation of the
multiplication law in $D(U)_q$. The new multiplication is to be a bilinear
map
$$*:\,D(U)_q \times D(U)_q \to D(U)_q[[t]],$$
$$*:\,f_1 \times f_2 \mapsto f_1 \cdot f_2+\sum_{i=1}^\infty
t^iC_i(f_1,f_2),$$
which satisfies the formal associativity condition
$$\sum_{i+k=m}C_i(f_1,C_k(f_2,f_3))=\sum_{i+k=m}C_i(C_k(f_1,f_2),f_3)$$
(cf. \cite{L}). When producing the new multiplication $*$, we follow F.
Berezin \cite{Be}. The bilinear operators $C_j:D(U)_q \times D(U)_q \to
D(U)_q$, $j \in{\Bbb N}$, will turn out to be $q$-bidifferential, and we
shall give explicit formulae for them.

 To begin with, choose a positive $\alpha$ and consider a linear functional
$\nu_ \alpha:{\rm Pol}({\Bbb C})_q \to {\Bbb C}$;
$$\int \limits_{U_q}fd \nu_ \alpha \stackrel{def}{=}{1-q^{4 \alpha}\over
1-q^2}\cdot \int \limits_{U_q}f \cdot (1-zz^*)^{2 \alpha+1}d \nu=(1-q^{4
\alpha})\,{\rm tr}\,T(f \cdot(1-zz^*)^{2 \alpha}).$$

Impose a norm $\Vert f \Vert_ \alpha=\left(\displaystyle \int
\limits_{U_q}f^*fd \nu_ \alpha \right)^{1/2}$ on ${\rm Pol}({\Bbb C})_q$.
Let $L_{q,\alpha}^2$ stand for the completion of ${\rm Pol}({\Bbb C})_q$
with respect to the above norm, and $H_{q,\alpha}^2$ for the linear span of
monomials $z^j \in L_{q,\alpha}^2$, $j \in{\Bbb Z}_+$.

\medskip

\begin{lemma}\hspace{-.5em}. The monomials $\{z^m\}_{m \in{\Bbb Z}_+}$ are
pairwise orthogonal in $H_{q,\alpha}^2$, and $\Vert z^m \Vert_
\alpha=((q^2;q^2)_m/(q^{4 \alpha+2};q^2)_m)^{1/2}$.
\end{lemma}

\smallskip

 {\bf Proof}. The pairwise orthogonality of the monomials $z^m$ is obvious;
$$\Vert z^m \Vert_ \alpha^2=(1-q^{4 \alpha})\cdot{\rm
tr}\,T(z^{*m}z^m(1-zz^*)^{2 \alpha})={1-q^{4 \alpha}\over
1-q^2}\int \limits_0^1(q^2y;q^2)_m \cdot y^{2 \alpha-1}d_{q^2}y=$$
$$={1-q^{4 \alpha}\over 1-q^2}\cdot{\Gamma_{q^2}(2 \alpha) \cdot
\Gamma_{q^2}(m+1) \over \Gamma_{q^2}(m+2 \alpha+1)}={(q^2;q^2)_m \over(q^{4
\alpha+2};q^2)_m}.$$

 We have used the well known \cite[\S 1.11]{GR} identity
$$\int \limits_0^1t^{\beta-1}\cdot(tq^2;q^2)_{\alpha-1}d_{q^2}t=
{\Gamma_{q^2}(\beta)\Gamma_{q^2}(\alpha)\over
\Gamma_{q^2}(\alpha+\beta)}.\eqno \Box$$

\medskip

\begin{corollary}\hspace{-.5em}. Let $\widehat{z}$ be the operator of
multiplication by $z$ in $H_{q,\alpha}^2$, and $\widehat{z}^*$ the conjugate
operator. Then $\widehat{z}$, $\widehat{z}^*$ are bounded, and
$$\widehat{z}^*\widehat{z}=q^2 \widehat{z}\widehat{z}^*+1-q^2+q^{4
\alpha}\cdot{1-q^2 \over 1-q^{4
\alpha}}(1-\widehat{z}\widehat{z}^*)(1-\widehat{z}^*\widehat{z}).\eqno(7.2)$$
\end{corollary}

\smallskip

 {\bf Proof} \ follows from
$$\widehat{z}:z^m \mapsto z^{m+1},\;m \in {\Bbb Z}_+;\quad \widehat{z}^*:1
\mapsto 0,\quad \widehat{z}^*:z^m \mapsto{1-q^{2m}\over 1-q^{4
\alpha+2m}}z^{m-1},\;m \in{\Bbb N}.$$

 In fact,
$$(1-\widehat{z}\widehat{z}^*)^{-1}:\,z^m \mapsto((q^{-2m}-q^{4
\alpha})/(1-q^{4 \alpha})z^m,$$
$$(1-\widehat{z}^*\widehat{z})^{-1}:\,z^m \mapsto((q^{-2m-2}-q^{4
\alpha})/(1-q^{4 \alpha})z^m.$$
Hence
$(1-\widehat{z}\widehat{z}^*)^{-1}=q^2(1-\widehat{z}^*\widehat{z})^{-1}-q^{4
\alpha}{\textstyle 1-q^2 \over \textstyle 1-q^{4 \alpha}}$.\hfill $\Box$

\medskip

 Lemma 7.1 and corollary 7.2 were proved in the work by S. Klimek and A.
Lesniewski \cite{KL} on two-parameter quantization of the disc beyond the
frameworks of perturbation theory.

 To every element $f=\sum a_{ij}z^iz^{*j}\in{\rm Pol}({\Bbb C})$ we
associate the linear operator $\widehat{f}=\sum
a_{ij}\widehat{z}^i\widehat{z}^{*j}$ in $H_{q,\alpha}^2$. The formal
deformation of the multiplication law in the algebra of functions in the
quantum disc will be derived via an application of "Berezin quantization
procedure" $f \mapsto \widehat{f}$ to the ordinary multiplication in the
algebra of linear operators.

 More exactly, (7.2) allows one to get a formal asymptotic expansion
$$\widehat{f_1}\cdot \widehat{f_2}=\widehat{f_1 \cdot
f_2}+\sum_{k=1}^\infty q^{4 \alpha k}{C_k \widehat{(f_1,f_2)}}.
\quad f_1,f_2
\in {\rm Pol}({\Bbb C})_q,$$
with $C_k:{\rm Pol}({\Bbb C})_q \times{\rm Pol}({\Bbb C})_q \to{\rm
Pol}({\Bbb C})_q$, $k \in{\Bbb N}$, bilinear maps. In this way, we get a
formal deformation
$$*:\,{\rm Pol}({\Bbb C})_q \times{\rm Pol}({\Bbb C})_q \to{\rm Pol}({\Bbb
C})_q[[t]];$$
$$f_1*f_2=f_1 \cdot f_2+\sum_{k=1}^\infty t^k \cdot C_k(f_1,f_2);\quad
f_1,f_2 \in{\rm Pol}({\Bbb C})_q.$$

 We present an explicit formula for the multiplication $*$, and thus also
for bilinear maps $C_k$, $k \in{\Bbb N}$. Let $\stackrel{\sim}{\Box}$
be a linear operator in ${\rm Pol}({\Bbb C})_q^{op}\otimes{\rm Pol}({\Bbb C})_q$
given by
$$\stackrel{\sim}{\Box}=q^{-2}(1-(1+q^{-2})z^*\otimes z+q^{-2}z^{*2}\otimes
z^2){\partial^{(r)}\over \partial z^*}\otimes{\partial^{(l)}\over \partial
z},$$
and $m:\,{\rm Pol}({\Bbb C})_q \times{\rm Pol}({\Bbb C})_q \to{\rm
Pol}({\Bbb C})_q$, $m:\,\psi_1 \otimes \psi_2 \to \psi_1 \psi_2$ the
multiplication in ${\rm Pol}({\Bbb C})_q$.

\medskip

\begin{theorem}\hspace{-.5em}. For all $f_1,f_2 \in{\rm Pol}({\Bbb C})_q$
$$f_1*f_2=(1-t)\sum_{j \in{\Bbb Z}_+}t^j \cdot
m(p_j(\stackrel{\sim}{\Box})f_1 \otimes f_2),$$
with
$$p_j(\stackrel{\sim}{\Box})=\sum_{k=0}^j{(q^{-2j};q^2)_k
\over(q^2;q^2)^2_k}q^{2k}\cdot \prod_{i=0}^{k-1}(1-q^{2i}((1-q^2)^2 \cdot
\stackrel{\sim}{\Box}+1+q^2)+q^{4i+2}).$$
\end{theorem}

\smallskip

 The proof can be found in \cite[Theorem 8.4]{SSV5}.

\medskip

 {\sc example 7.4}. For all $f_1,f_2 \in{\rm Pol}({\Bbb C})_q$
$$f_1*f_2=f_2 \cdot f_2+t \cdot(q^{-2}-1){\partial^{(r)}f_1 \over \partial
z^*}(1-z^*z)^2{\partial^{(l)}f_2 \over \partial z}+O(t^2).$$

\medskip

 It is worthwhile to note that the formal associativity of the
multiplication $*$ follows from the associativity of multiplication in the
algebra of linear operators.

\medskip \stepcounter{theorem}

\begin{corollary}\hspace{-.5em}. The bilinear operators $C_k$ are of the
form (7.1) and are extendable by continuity up to $q$-bidifferential
operators $C_k:\,D(U)_q \times D(U)_q \to D(U)_q$.
\end{corollary}

\medskip

 The above $q$-bidifferential operators determine a formal deformation of
the multiplication in $D(U)_q$. The formal associativity of the newly formed
multiplication $*:\,D(U)_q \times D(U)_q \to D(U)_q[[t]]$ follows from the
formal associativity of the previous multiplication $*:\,{\rm Pol}({\Bbb
C})_q \times {\rm Pol}({\Bbb C})_q \to{\rm Pol}({\Bbb C})_q[[t]]$.

 Finally, note that our proof of theorem 7.3 is based on the properties of
some $q$-analogue for Berezin transformation \cite{UU}.

\bigskip

\section*{Appendix. On $q$-analogue of the Green formula}

 Consider the two-sided ideal $J \in{\rm Pol}({\Bbb C})_q$ generated by the
element $1-zz^*\in{\rm Pol}({\Bbb C})_q$, and the commutative quotient
algebra ${\Bbb C}[\partial U]_q \stackrel{def}{=} {\rm Pol}({\Bbb C})_q/J$.
Its elements will be identified with the corresponding polynomials on the
circle $\partial U$. The image $f|_{\partial U}$ of $f \in{\rm Pol}({\Bbb
C})_q$ under the canonical homomorphism ${\rm Pol}({\Bbb C})_q \to{\Bbb
C}[\partial U]_q$ will be called a restriction of $f$ onto the boundary of
the quantum disc.

 Define the integral $\Omega({\Bbb C})_q^{(1,0)}\to{\Bbb C}$ by
$$\int \limits_{\partial U}dz \cdot f \stackrel{def}{=}2 \pi i \int
\limits_{\partial U}(z \cdot f)|_{\partial U}d \nu,\quad f \in {\rm
Pol}({\Bbb C})_q,$$
with
$$\int \limits_{\partial U}\psi d \nu \stackrel{def}{=}\int \limits_0^{2
\pi}\psi(e^{i \theta}){d \theta \over 2 \pi}.$$

\medskip

{\bf Proposition A.1.} {\it For all $\psi \in \Omega({\Bbb C})_q^{(0,1)}$
one has
$$\int \limits_{U_q}\overline{\partial}\psi=\int \limits_{\partial
U}\psi.\eqno(A.1)$$}

\smallskip

 {\sc Remark A.2.} The integral $\int \limits_{U_q}dz \cdot f \cdot dz^*=-2i \pi
\int \limits_{U_q}fd \mu$ introduced in section 4 for $f \in D(U)_q$, is
extendable by continuity onto all (1,1)-forms $dz \cdot f \cdot dz^*$ with
$$f=\sum_{m>0}z^m
\psi_m(y)+\psi_0(y)+\sum_{m>0}\psi_{-m}(y)z^{*m}\in D(U)_q^\prime,$$
such that $\sum \limits_{m \in{\Bbb Z}_+}|\psi_0(q^{2m})|q^{2m}<\infty$. Under
these assumptions one also has
$$\int \limits_{U_q}dzdz^*f=\int \limits_{U_q}dz \cdot f \cdot
dz^*=\int \limits_{U_q}fdzdz^*.\eqno(A.3)$$

\smallskip

 {\bf Proof}. Recall that (see (1.2))
$$\psi=dz \left(\sum_{m>0}z^m
\psi_m(y)+\psi_0(y)+\sum_{m>0}\psi_{-m}(y)z^{*m}\right).\eqno(A.2)$$
We can restrict ourselves to the case $\psi=dz \psi_{-1}(y)z^*$, since this
is the only term in (A.2) which could make a non-zero contribution to (A.1).

 It follows from the definitions that
$\overline{\partial}\psi=dz \cdot f(y)dz^*$, with
$$f(y)=\psi_{-1}(y)-q^{-2}
\frac{\psi_{-1}(q^{-2}y)-\psi_{-1}(y)}{q^{-2}y-y}(1-y).$$
In fact, $\overline{\partial}y=\overline{\partial}(1-zz^*)=-zdz^*$. Hence
$\overline{\partial}y^m=\displaystyle \sum
\limits_{j=0}^{m-1}y^j(-zdz^*)y^{m-1-j}=-{\textstyle 1-q^{2m}\over
\textstyle 1-q^2}zy^{m-1}dz^*$.

 That is, for any polynomial $p(y)$ one has
$$\overline{\partial}p(y)=-z{p(y)-p(q^2y) \over y-q^2y}\cdot
dz^*.\eqno(A.4)$$
(Note that the validity of (A.4) for polynomials already implies its
validity for all distributions). Finally,
$$\overline{\partial}(dz
\psi_{-1}(y)z^*)=dz \left(-z{\psi_{-1}(y)-\psi_{-1}(q^2y)\over
y-q^2y}z^*+\psi_{-1}(y)\right)dz^*.$$
On the other hand, $-z{\textstyle \psi_{-1}(y)-\psi_{-1}(q^2y)\over
\textstyle y-q^2y}z^*+\psi_{-1}(y)=f(y)$, since $zy=q^{-2}yz$, $zz^*=1-y$.

 If one assumes $\psi_{-1}(0)=0$, it is easy to show that $\sum
\limits_{n \in{\Bbb Z}_+}f(q^{2n})q^{2n}=0$. Hence, in this case $\int
\limits_{U_q}\overline{\partial}\psi=\int \limits_{\partial U}\psi=0$. Thus,
Proposition A.1 is proved for all (1,0)-forms from some linear subspace of
codimensionality 1. Now it remains to prove (A.1) in the special case
$\psi=dz \cdot z^*$. \hfill $\Box$

\medskip

 {\bf Corollary A.3}. {\it If $\psi \in \Omega(U)_q^{(1,0)}$, then $\int
\limits_{U_q}\overline{\partial}\psi=0$}.

\end{document}